\documentclass[11pt]{amsart}
\usepackage{amscd}
\usepackage{hyperref}
\def\and{\operatorname{and}}
\def\ker{\operatorname{ker}}

\def\dim{\operatorname{dim}}
\def\depth{\operatorname{depth}}
\def\grade{\operatorname{grade}}
\def\ass{\operatorname{Ass}}
\def\modulo{\operatorname{modulo}}

\newcommand{\m}{\frak m}
\newcommand{\ZZ}{\mathbb Z}

\newcommand{\R}{\mathcal R}

\newcommand{\F}{\mathcal F}

\newtheorem{lemma}{Lemma}[section]
\newtheorem{corollary}[lemma]{Corollary}
\newtheorem{theorem}[lemma]{Theorem}
\newtheorem{proposition}[lemma]{Proposition}

\newtheorem{definition}[lemma]{Definition}
\newtheorem{remark}[lemma]{Remark}

\newtheorem{example}[lemma]{Example}

\advance\headheight1.15pt
\begin{document}
\title{Hilbert coefficients and depth of fiber cones}

\author{A. V. Jayanthan}
\address{School of Mathematics, Tata Institute of Fundamental
Research,
Homi Bhabha Road, Colaba, Mumbai, India - 400005}
\thanks{The
first author was supported by the National Board for Higher Mathematics, India}
\thanks{AMS Classification 2000: 13H10, 13H15 (Primary), 13C15, 13A02
(Secondary)}
\email{jayan@math.tifr.res.in}
\author{J. K. Verma}
\address{Department of Mathematics, Indian Institute of Technology
Bombay, Powai, Mumbai, India - 400076}
\email{jkv@math.iitb.ac.in}

\maketitle
\centerline{
{\em Dedicated to Prof. Wolmer Vasconcelos on the occasion of his 65th 
birthday}}

\begin{abstract}
Criteria are given in terms of certain Hilbert coefficients for the
fiber cone $F(I)$ of an $\m$-primary ideal $I$ in a Cohen-Macaulay
local ring $(R,\m)$ so that it is Cohen-Macaulay or has depth at least
$dim(R)-1.$ A version of Huneke's fundamental lemma is proved for
fiber cones. Goto's results concerning Cohen-Macaulay fiber cones
of ideals with minimal multiplicity are obtained  as consequences. 
\end{abstract}

\section{Introduction}
Let $(R,\m)$ be a $d$-dimensional Cohen-Macaulay local ring  
having infinite residue field $R/\m.$ Let $I$ be 
an $\m$-primary ideal. The objective of this paper 
is to explore some connections between  the 
coefficients of the polynomial $P_{\m}(I,n)$ 
corresponding to the function 
$H_{\m}(I,n)=\lambda(R/ \m I^n)$ and depth of the 
fiber cone $F(I)=\oplus_{n \geq 0}I^n/\m I^n$ of $I$.
The relation between Hilbert coefficients and depth  
 has been a subject of several papers in 
the context of associated graded rings and Rees algebras 
of ideals. That  conditions on Hilbert coefficients could 
force high depth for associated graded rings was first 
observed by Sally in \cite {s1}. Since then
numerous conditions have been provided for  the Hilbert
coefficients so that the associated graded ring of $I$, 
$G(I)=\oplus_{n \geq 0}I^n/I^{n+1}$, is either Cohen-Macaulay 
or has almost maximal depth, i.e. the grade of the maximal 
homogeneous ideal of $G(I)$ is at least $d-1.$ 

Let $J$ be a minimal reduction of $I.$ In their elegant 
paper \cite{hm} Huckaba and Marley provided
necessary and sufficient conditions on the coefficients
of the  Hilbert polynomial $P(I,n)$ corresponding to the 
Hilbert function $H(I,n)=\lambda(R/I^n)$ so that $G(I)$ is
Cohen-Macaulay and has almost maximal depth. We shall write the
Hilbert polynomial $P(I,n)$ in the following way:

$$
P(I,n)=e_0(I){ n+d-1 \choose d}-e_1(I){ n+d-2 \choose d-1}
       + \cdots + (-1)^d e_d.
$$

Huckaba and Marley  showed:\\   
(i) $G(I)$ is Cohen-Macaulay if and only if 
$e_1(I)=\sum_{n \geq 1 }\lambda(I^{n}+J/J)$ 
 \\(i) $G(I)$ has almost  maximal depth if 
and only if $e_1(I)=\sum_{n \geq 1 }\lambda(I^{n}/JI^{n-1}).$ 

Their results unify several results known in the literature and
provide an effective approach to deal with such questions. 

The relation between Hilbert coefficients and depth of $F(I)$ is not
well understood. The papers \cite{cz}, \cite{g}, \cite{dv} and
\cite{drv} provide sufficient conditions in terms of certain Hilbert
coefficients for the Cohen-Macaulay property of $F(I).$ Recently, a
number of papers have appeared which discuss Cohen-Macaulay and other
properties of fiber cones. See for example \cite{cgpu}, \cite{c},
\cite{cpv}, \cite{hh}.  The form ring $G(I)$ and the fiber cone $F(I)$
coincide when $I=\m.$ This indicates that there may be  analogues of
results of Huckaba and Marley for the fiber cone. The first guess for
the appropriate Hilbert function to be used is naturally the Hilbert
function of $F(I)$.  However we have observed that this does not seem
to be of much help in predicting depth. We will show that the
coefficients of the polynomial $P_{\m}(I,n)$  corresponding to the
function  $\lambda(R/\m I^n)$ control the Cohen-Macaulay and almost
maximal depth properties of $F(I).$     
   
We now describe the main results of this paper. Write the polynomial 
$P_{\m}(I,n)$ as

$$ 
P_{\m}(I,n) =g_0(I){ n+d-1 \choose d}-g_1(I){ n+d-2 \choose d-1}
       + \cdots + (-1)^d g_d(I).
$$

\noindent
Let grade $G(I)_+ \geq d-1. $ In sections $4$ and $5$ 
we shall prove that \\ (i) $F(I)$ is Cohen-Macaulay if and only  if 
$g_1(I)= \sum_{n \geq 1} \lambda(\m I^n + JI^{n-1}/JI^{n-1}) -1$
and \\ (ii) $F(I)$ has almost maximal depth if and only if 
$g_1(I)= \sum_{n \geq 1} \lambda(\m I^n / \m JI^{n-1}) -1.$

It can be seen that the minimal number of generators of $I$, $\mu(I)
\leq e_0(I)+d - \lambda(R/I). $ Goto in \cite{g} defined an ideal
$I$ in a Cohen-Macaulay local ring to have {\em minimal multiplicity}
if equality holds in the above inequality. He showed that if $I$ has
minimal multiplicity then $F(I)$ is Cohen-Macaulay if and only if
$G(I)$ has almost maximal depth. We shall recover this result in
section $6$ as a consequence of our criterion for Cohen-Macaulayness
in terms of $g_1(I).$ In fact we shall prove that $I$ has minimal
multiplicity if and only if $g_1(I)=-1.$ We will also recover
a bound for the multiplicity of the fiber cone due to Corso, Polini
and Vasconcelos \cite{cpv}.
 
Since the criteria for Cohen-Macaulay and almost maximal depth
properties of $F(I)$ require one to know the coefficient $g_1(I)$, it
is desirable to have an effective method of its computation. In
Section $5$ we shall show that in a one dimensional Cohen-Macaulay
local ring $g_1(I)=\sum_{n \geq 1}\lambda(\m I^n/J\m I^{n-1})-1.$ We
will also present a simple proof of a criterion due to  Cortadellas
and Zarzuela \cite{cz} for the sequence of initial forms in $F(I)$
of elements in a regular sequence in $R$ to be a regular sequence in
$F(I).$ 

In Section $3$ 
we shall generalize  the  {\em fundamental lemma of Huneke } in \cite{h2} for  
finding a formula for  $g_1(I).$ However, we need a modified version 
of this lemma so that it works for the function $\lambda(R/ \m I^n).$

In the second section we will discuss the technical topic of superficial 
and regular elements in fiber cones.

As no extra effort is required, we will  prove our results for filtrations of 
ideals. In a subsequent paper on fiber cones we will see that it is useful 
to develop the 
criteria for filtrations as sometime we need to deal with filtrations to 
prove results about the $I$-adic filtration.

{\em Acknowledgements:} We thank Balwant Singh and Santiago Zarzuela 
for discussions.

\section{Superficial and regular elements in fiber cones of
filtrations}
In this section we will gather some results on superficial and 
regular elements in fiber cones.
Throughout this paper $(R,\m)$ will denote a Noetherian 
local ring of positive (Krull) dimension $d$, with maximal ideal 
$\m$  and infinite residue field $R/\m$. 
A sequence of ideals  $\F = \{I_n\}_{n\geq 0}$ 
is called a filtration 
if $I_0=R$, $I_1 \neq R$, $I_{n+1} \subseteq I_n$, 
and $I_n I_m \subseteq I_{n+m}$ for all $n, m \geq 0$. 
The Rees algebra $\R(\F)$ and the 
associated graded ring $G(\F)$ are defined to be the graded rings:
$$
\R(\F)=R \oplus I_1t\oplus I_2t^2\oplus \cdots , 
     G(\F)=R/I_1\oplus I_1/I_2 \oplus I_2/I_3 \oplus \cdots
$$

The filtration $\F$ is called Noetherian if $\R(\F)$ is a 
Noetherian ring. Throughout the paper we will assume that $\F$ is 
Noetherian and $I_n \neq 0$ for all $n \geq 0.$ The ideal generated by 
elements of positive degree in $G(\F)$ will be denoted by $G(\F)_+$.
The filtration $\F$ is called $I_1$-good if $\R(\F)$ is a finite module 
over the Rees algebra $\R(I_1).$ An $I_1$-good filtration is called 
a Hilbert filtration if $I_1$ is $\m$-primary.

The fiber cone of $\F$ with respect to an ideal $K$ containing $I_1$
is the graded ring 
$$F_K(\F) = R/K \oplus I_1/KI_1 \oplus I_2/KI_2 \oplus \cdots .$$
  
For $x \in I_1$, let $x^*$ and $x^o$ denote the
initial form in degree one component of $G(\F)$ and $F_K(\F)$
respectively. We will always assume that $I_{n+1}\subseteq KI_n $ for all
$n \geq 0.$ This is required in Lemma \ref{supchar} which is essential
in all the arguments that apply induction on the dimension of $R$ in 
the subsequent sections. 

\begin{definition}
For an element $x \in I_1$ such that $x^o \neq 0 \in F_K(F)$, $x^o$ is
said to be superficial in $F_K(\F)$ if there
exists $c > 0$ such that $(0 : x^o) \cap F_K(\F)_{n} = 0$ for all 
$n > c.$
\end{definition}

\noindent
It can easily be seen that $x^o$ is superficial in $F_K(\F)$ if and
only if there exists $c > 0$ such that 
$(KI_{n+1} : x) \cap I_{n} = KI_{n}$ for all $n > c$.  
We first show the existence of superficial
elements in $F_K(\F)$ and proceed to prove some of their properties.
The existence of superficial elements in a graded ring is well-known.
Since in our case we need existence of elements which are superficial
in $F_K(\F)$ as well as $G(\F)$ simultaneously, we give a proof of the
following result for the sake of completeness. 

\begin{proposition}
Let $(R,\m)$ be a Noetherian local ring of dimension $d > 0$. Let $\F$
be an $I_1$-good filtration and $K$ be an $\m$-primary ideal
containing $I_1$. Then there exists $x \in I_1 \backslash KI_1$ such
that $x^o$ is superficial in $F_K(\F)$ and $x^*$ is superficial in
$G(\F)$.
\end{proposition}
\begin{proof}
Let the set of associated primes of $G(\F)$ and $F_K(\F)$ be 
$$
\ass(G(\F)) = \{P_1, \ldots, P_r, P_{r+1}, \ldots, P_s\} \and 
\ass(F_K(\F)) = \{Q_1, \ldots, Q_l, Q_{l+1}, \ldots, Q_m \}
$$
such that for all $n \gg 0$, $I_n/I_{n+1} \subseteq P_i$ for 
$r+1 \leq i \leq s$ and $I_n/KI_n \subseteq Q_j$ for 
$l+1 \leq j \leq m$. The associated graded ring $G(\F)$ 
and the fiber cone $F_{K}(\F)$ are 
both homomorphic images of the extended Rees algebra 
$\R(\F)(t^{-1})$  since  $G(\F)=\R(\F)(t^{-1})/(t^{-1})$ and   
$F_{K}(\F)=  \R(\F)(t^{-1})/(t^{-1}, K).$  
Let $\mathcal {P} = \{P_1', \ldots, P_s', Q_1', \ldots, Q_m'\}$ be the
collection of prime ideals in the extended Rees algebra $\R(\F)(t^{-1})$
 which are the pre-images of the
ideals $P_1, \ldots, P_s $  in $\ass(G(\F))$ and  
$Q_1, \ldots, Q_m  $ in  $\ass(F_K(\F))$ respectively.  
Since $ R/\m $ is infinite, 
$\R_1 \neq \m I_1t \cup_{i=1}^r P_i' \cup_{i = 1}^l Q_i'$.
Choose $xt \in \R_1 \backslash \m I_1t \cup_{i=1}^r P_i' \cup_{i =
1}^l Q_i'$.  We show that $0 \neq x^o \in F_K(\F)_1$ is superficial in
$F_K(\F)$ and $x^* \in G(\F)_1$ is superficial in $G(\F)$.  We need to
show that $(0 : x^o) \cap F_K(\F)_n = 0$ for $n \gg 0$.  Let $y^o \in
(0 : x^o)$. Let $(0) = N_1 \cap \cdots \cap N_m$ be a primary
decomposition of $(0)$ in $F_K(\F)$ such that $N_i$ is $Q_i$-primary
for $i = 1, \ldots, m$. Then $y^ox^o \in N_i$ for all $1 \leq i \leq
l$ and $x^o \notin Q_i$ for $i = 1, \ldots, l$. Therefore $y^o \in
N_i$ for $i = 1, \ldots, l$. Thus $(0 : x^o) \subseteq N_1 \cap \cdots
\cap N_l$. For $l+1 \leq j \leq m$, $F_K(\F)_n \subseteq Q_j$ for $n
\gg 0$. Therefore there exists a $c > 0$ such that $\oplus_{n\geq
c}F_K(\F)_n \subseteq N_{l+1} \cap \cdots \cap N_m$.  Therefore for
all $n \geq c$
$$
(0 : x^o) \cap F_K(\F)_n \subseteq N_1 \cap \cdots \cap N_l \cap 
N_{l+1} \cap \cdots \cap N_m  = (0).
$$
Therefore $x^o$ is superficial in $F_K(\F)$. A similar argument shows
that $x^*$ is superficial in $G(\F)$.
\end{proof}
\noindent
In the next lemma, we characterize the
property of  an element being superficial in the fiber cone in terms of 
its properties in the local ring.
 
\begin{lemma}\label{supchar}
(i) If there exists a $c > 0$ such that $(KI_n : x) \cap I_c =
KI_{n-1}$ for all $n > c$, then $x^o$ is superficial in $F_K(\F)$. 
\vskip 2mm
\noindent
(ii) If $x^o$ is superficial in $F_K(\F)$ and $x^*$ is
superficial in $G(\F)$, then there exists $c > 0$ such that $(KI_n :
x) \cap I_c = KI_{n-1}$ for all $n > c$.  Moreover if $x$ is regular
in $R$, then $KI_n : x = KI_{n-1}$ for all $n \gg 0$. 
\end{lemma}
\begin{proof}
$(i)$ Suppose $(KI_{n} : x) \cap I_c = KI_{n-1}$ for all $n > c$.
Then $(KI_{n} : x) \cap I_{n-1} \subseteq (KI_{n} : x) \cap I_c =
KI_{n-1}$ for all $n > c$. Therefore $x^o$ is superficial in
$F_K(\F)$. 
\vskip 2mm
\noindent
$(ii)$ Suppose $x \in I_1$ is such that $x^o$ is superficial in
$F_K(\F)$ and $x^*$ is superficial in $G(\F)$. Then there exist $c_1,
c_2$
such that for all $n > c_1$, $(0 : x^*) \cap G(\F)_n = 0$ and $(0 : x^o)
\cap F_K(\F)_n = 0$ for all $n > c_2$. Choose $c = \max\{c_1,
c_2\}+1$. 
\vskip 2mm
\noindent
{\it Claim :} $(KI_n : x) \cap I_c = KI_{n-1}$ for all $n > c$. 
\vskip 2mm
\noindent
Suppose that there exists an integer $t$ such that $y \in I_t 
\backslash KI_t$. If $xy \in KI_{t+1}$, then $y^o \in (0
: x^o) \cap F_K(\F)_t = 0$ which is a contradiction (the last equality
holds because $t \geq c > c_2$). Therefore, $xy \notin
KI_{t+1}$. Since $xy \in KI_n$, $n < t+1$ so that $n \leq t$.
Therefore $y \in I_t \subseteq I_n \subseteq KI_{n-1}$. 
\vskip 2mm
\noindent
Now, suppose that the image of $y$ in $F_K(\F)$ is zero, i.e., given
any $n$ such that $y \in I_n$, $y \in KI_n$. Then $y \in (KI_n : x)
\cap I_c \subseteq (I_n : x) \cap I_c = I_{n-1}$, since $c > c_1$. 
Therefore, $y \in KI_{n-1}$.
\vskip 2mm
\noindent
Suppose that $x$ is regular in $R$. Then by the Artin-Rees lemma,
there exists a $c$ such that $KI_n \cap (x) = I_{n-c}(KI_c \cap (x))
\subseteq xI_{n-c} \subseteq xI_c$ for all $n \geq 2c$. 
Therefore $KI_n : x \subseteq I_c$
for $n \gg 0$. Hence for all $n \gg 0$, $KI_n : x = (KI_n : x) \cap
I_c = KI_{n-1}$.
\end{proof}

\noindent
For an element $x \in I$ such that $x^*$ is superficial in $G(I)$, it
is known that $x \in I \backslash I^2$. In the following result we
show that a similar property is true for superficial elements in
$F_K(\F)$.
\begin{lemma}
If $x^o \in I_1/KI_1$ is superficial in $F_K(\F)$ and $x^*$ is
superficial in $G(\F)$, then $x \in I_1 ~ \backslash~  KI_1$.
\end{lemma}
\begin{proof}
Since $x^o \in F_K(\F)$ and $x^* \in G(\F)$ are superficial, by Lemma
$\ref{supchar}$, there exists $c > 0$ such that $(KI_n : x) \cap I_c =
KI_{n-1}$ for all $n > c$. Put $n = c+1$. Then $(KI_{c+1} : x) \cap
I_c = KI_c$. Suppose $x \in KI_1$.  Let $y \in I_c$. Then $xy \in
KI_{c+1}$ so that $y \in (KI_{c+1} : x) \cap I_c = KI_c$. Therefore
$I_c = KI_c$. By Nakayama Lemma $I_c = 0$ which is a contradiction to
the fact that $\F$ is a Hilbert filtration.  Therefore $x \in I_1 ~
\backslash~ KI_1$.
\end{proof}

For the fiber cone $F_K(\F)$, let $H(F_K(\F), n) = \lambda(F_K(\F)_n)
= \lambda (I_n/KI_n)$ denote its Hilbert function and let $P(F_K(\F),
n)$ denote the corresponding Hilbert polynomial.
\begin{proposition}
Let $x^o \in F_K(\F)$ be superficial. Then 
$\dim F_K(\F)/x^oF_K(\F) = \dim F_K(\F) - 1$.
\end{proposition}
\begin{proof}
Consider the exact sequence 
$$
0 \longrightarrow (KI_n : x) \cap I_{n-1}/KI_{n-1}
\longrightarrow I_{n-1}/KI_{n-1} \stackrel{x}{\longrightarrow}
I_n/KI_n \longrightarrow I_n/(KI_n + xI_{n-1})
\longrightarrow 0.
$$
Then $H(F_K(\F),n) - H(F_K(\F), n-1) = H(F_K(\F)/x^oF_K(\F), n) -
\lambda((KI_n : x) \cap I_{n-1}/KI_{n-1})$ for all $n \geq 1$. Since
$x^o$ is superficial in $F_K(\F)$, $(KI_n : x) \cap I_{n-1} =
KI_{n-1}$ for $n \gg 0$, so that $P(F_K(\F)/(x^o), n) =
P(F_K(\F), n) - P(F_K(\F), n-1)$. Hence $\dim F_K(\F)/x^oF_K(\F) =
\dim F_K(\F) - 1$.  
\end{proof}

\noindent
In the following lemma we provide a characterization for regular
elements in $F_K(\F)$. It can be seen that this property is quite
similar to the behaviour of regular elements in $G(I)$.
\begin{lemma}
For $x \in I_1 \backslash KI_1$, $x^o \in F_K(\F)$ is regular if and
only if $(KI_n : x) \cap I_{n-1} = KI_{n-1}$ for all $n \geq 1$.  If
$x^*$ is regular in $G(\F)$ and $x^o$ is regular in $F_K(\F)$ then
$KI_n : x = KI_{n-1}$ for all $n \geq 1.$
\end{lemma}
\begin{proof}
Suppose $(KI_n : x) \cap I_{n-1} = KI_{n-1}$ for all $n \geq 1$. In
other words, $(0 : x^o) \cap F_K(\F)_{n-1} = (0)$ for all $n \geq 1$.
Since $(0 : x^o)$ is a homogeneous ideal, $(0 : x^o) = \oplus_{n \geq
0} (0 : x^o) \cap F_K(\F)_n = (0)$. Conversely, assume that $x^o$ is a
regular element in $F_K(\F)$. Then $\alpha_n : F_K(\F)_{n-1}
\longrightarrow F_K(\F)_n$ is an injective map for all $n \geq 1$,
where $\alpha_n$ is the multiplication by $x^o$ for all $n \geq 1$.
Since $\ker \alpha_n = (KI_n : x) \cap I_{n-1}/KI_{n-1} = 0$, $(KI_n :
x) \cap I_{n-1} = KI_{n-1}$ for all $n \geq 1$.
\vskip 2mm
\noindent
If $KI_n : x = KI_{n-1}$ for all $n \geq 1$, then clearly $x^o$ is
regular in $F_K(\F)$. Suppose that $x^*$  is regular in $G(\F)$ and
$x^o$ is regular in $F_K(\F)$. Let $y \in KI_n : x$. If there exists a
$t$ such that $0 \neq y^o \in I_t/KI_t$, then $0 \neq y^ox^o \in
I_{t+1}/KI_{t+1}$ so that $yx \notin KI_{t+1}$. Since $yx \in KI_n$,
$n < t+1$. Therefore $y \in I_t \subseteq I_{n} \subseteq KI_{n-1}.$
Suppose we can not find $t$ such that $0 \neq y^o \in I_t/KI_t$. Then,
if $y \in I_n$, $y \in KI_n$. Since $y \in KI_n : x$, $y \in I_n : x$.
By hypothesis, $x^*$ is regular in $G$. Hence $I_n : x = I_{n-1}$ for
$n \geq 1$.  Therefore $y \in I_{n-1}$ and hence $y \in KI_{n-1}$.
\end{proof}

\noindent
The following lemma is an analogue of Lemma 2.2 of \cite{hm}.
This will play a crucial role in induction arguments. This is the so-called
Sally-machine for fiber cones.

\begin{lemma}\label{sally}
Let $x \in I_1$ be such that $x^*$ is superficial in $G(\F)$ and $x^o
\in F_K(\F)$ is superficial in $F_K(\F)$. Let $\bar{\F} = \{I_n +
xR/xR\}_{n\geq 0}$ and $\bar{K} = K/xR$. If $\depth
F_{\bar{K}}(\bar{\F}) > 0$, then $x^o$ is regular in $F_K(\F)$.
\end{lemma}
\begin{proof}
Let $y^o \in I_t/KI_t$ be such that its natural image $\bar{y}^o$ is a
regular element in $F_{\bar{K}}(\bar{\F})$. Then $(KI_{n+tj} : y^j)
\cap I_n \subseteq (KI_n, x)$ for all $n, j \geq 1$. Since $x^o$ is
superficial in $F_K(\F)$ and $x^*$ is superficial in $G(\F)$, there 
exists $c > 0$ such that $(KI_{n+j} : x^j) \cap I_c = KI_n$ for all 
$n > c \mbox{ and } j \geq 1$, by Lemma \ref{supchar} (ii).  
Let $n$ and $j$ be arbitrary and $p > c/t$. Then $y^p(KI_{n+j} : 
x^j) \subseteq
(KI_{n+pt+j} : x^j) \cap I_c \subseteq KI_{n+pt}$. Thus 
$$
(KI_{n+j} : x^j) \cap I_n \subseteq (KI_{n+pt} : y^p) \cap I_n
\subseteq (KI_n, x).
$$ 
Therefore $(KI_{n+j} : x^j) \cap I_n \subseteq KI_n + x(KI_{n+j}
: x^{j+1})$. Iterating this formula $n+1$ times we get,
\begin{eqnarray*}
(KI_{n+j} : x^j) \cap I_n & \subseteq & KI_n + xKI_{n-1} + \cdots + 
x^{n+1}(KI_{n+j} : x^{n+j+1}) \\
 & = & KI_n.
\end{eqnarray*}
Therefore $(KI_{n+j} : x^j) \cap I_n = KI_n$ for all $n \geq 1$ and 
hence $x^o$ is regular in $F_K(\F)$.
\end{proof}
\noindent
\section{Hilbert coefficients for the function $\lambda(R/KI_n)$}

Throughout this section  $\F = \{I_n\}_{n\geq
0}$ will be a Hilbert filtration of $R$. Let $K$ be an ideal of $R$
such that $I_{n+1} \subseteq KI_n$ for all $n \geq 0$. 
Let $H(F, n) = \lambda(F_K(\F)_n) = \lambda(I_n/KI_n)$ be
the Hilbert function of the fiber cone $F = F_K(\F)$. Then, $H(F, n) =
\lambda (R/KI_n) - \lambda (R/I_n)$. Since both $H(F,n)$ and 
$\lambda(R/I_n)$ are polynomials for $n \gg 0$, $\lambda(R/KI_n)$ is
also a polynomial for $n \gg 0$. Since the coefficients of this
polynomial are related with the Hilbert coefficients of the fiber cone
and the Hilbert-Samuel coefficients of $\F$, it is expected that their
properties will be related with the properties of the fiber cone.
Huneke's fundamental lemma \cite{h2} provides formulas for the Hilbert
coefficients of the Hilbert polynomial of an $\m$-primary ideal in a
two-dimensional Cohen-Macaulay local ring. We will prove an analogue
of this lemma for the fiber cones. It will yield formulas for the
Hilbert polynomial of the fiber cone once we have access to a minimal
reduction of $I_1.$

\vskip 4mm
\noindent
Let $H_K(\F, n) = \lambda(R/KI_n)$ (resp. $H(\F, n) = \lambda(R/I_n)$)
and let $P_K(\F,n)$ (resp. $P(\F, n)$) be the corresponding
polynomial. Since $P_K(\F,n) = P(\F,n) + P(F,n)$, it is a polynomial
of degree $d$ with leading coefficient $e_0(I)$. We write the above
polynomials in the following way:
\begin{eqnarray*}
P(F,n) & = & f_0 {n+d-2 \choose d-1} - f_1 {n+d-3 \choose d-2} + 
\cdots + (-1)^{d-1}f_{d-1}, \\
P(\F,n) & = & e_0 {n+d-1 \choose d} - e_1 {n+d-2 \choose d-1} + 
\cdots + (-1)^{d}e_{d}, \\
P_K(\F,n) & = & g_0 {n+d-1 \choose d} - g_1 {n+d-2 \choose d-1} + 
\cdots + (-1)^{d}g_{d}.  
\end{eqnarray*}

Then $g_0 = e_0$ and $g_i = e_i - f_{i-1}$ for all $1 \leq i \leq d$.
Note that $P(F,n)$, given above, is not written in the standard way of
representing the Hilbert function of a graded algebra. The changed
notation is for the convenience of the computation. 
For a numerical function $h : \ZZ \longrightarrow \ZZ$, let $\Delta
h(n) := h(n) - h(n-1)$.

\begin{lemma}\label{funda}
Let $(R,\m)$ be a $2$-dimensional Cohen-Macaulay local ring. Let $\F$
be a Hilbert filtration, $K$ be an ideal with $I_1 \subseteq K$ and $J
= (x, y)$ be a minimal reduction of $I_1$. Then for all $n \geq 2$,
$$
\Delta^2\left[P_K(\F,n) - H_K(\F,n)\right] =
\lambda\left(\frac{KI_n}{KJI_{n-1}}\right) - 
\lambda\left(\frac{KI_{n-1} : J}{KI_{n-2}}\right).
$$
\end{lemma}
\begin{proof}
Consider the exact sequence
$$
\begin{CD}
0 @>>> \frac{R}{KI_{n-1} : J} @>\beta>> \left(\frac{R}{KI_{n-1}}\right)^2
@>\alpha>> \frac{J}{KJI_{n-1}} @>>> 0,
\end{CD}
$$
where $\alpha$ is the map $\alpha(\bar{r}, \bar{s}) = 
\overline{xr + ys}$ and $\beta(\bar{r}) = (\bar{y}\bar{r},
-\bar{x}\bar{r})$.
It follows that for all $n \geq 2$
\begin{eqnarray*}
2\lambda(R/KI_{n-1}) & = & \lambda(R/(KI_{n-1} : J)) +
\lambda(J/KJI_{n-1}) \\
& = & \lambda(R/(KI_{n-1} : J)) + \lambda(R/KJI_{n-1}) - \lambda(R/J).
\end{eqnarray*}
Therefore $e_0(\F) + 2\lambda(R/KI_{n-1}) = \lambda(R/KJI_{n-1}) +
\lambda(R/(KI_{n-1} : J))$. Hence
\begin{eqnarray*}
e_0(\F) & - & \lambda(R/KI_{n}) + 2\lambda(R/KI_{n-1}) - 
\lambda(R/KI_{n-2})\\
& = & \lambda(R/KJI_{n-1}) - \lambda(R/KI_{n}) + \lambda(R/(KI_{n-1} :
J))- \lambda(R/KI_{n-2}) \\
& = & \lambda(KI_n/KJI_{n-1}) - \lambda(KI_{n-1} : J/KI_{n-2})
\end{eqnarray*}
Since $\Delta^2P_K(\F,n) = e_0(\F)$, 
$$
\Delta^2\left[P_K(\F,n) - H_K(\F,n)\right] =
\lambda\left(\frac{KI_n}{KJI_{n-1}}\right) - 
\lambda\left(\frac{KI_{n-1} : J}{KI_{n-2}}\right).
$$
\end{proof}
\begin{corollary}[\cite{h2}, Fundamental Lemma 2.4]
Let $(R,\m)$ be a $2$-dimensional Cohen-Macaulay local ring and $I$ be
an $\m$-primary ideal. Let $J = (x,y)$ be a minimal reduction of $I$.
Let $H(I,n) = \lambda(R/I^n)$ be the Hilbert function of $I$ and let
$P(I,n)$ be the corresponding Hilbert polynomial. Then for all $n \geq
2$,
$$
\Delta^2[P(I,n) - H(I,n)] = \lambda\left(\frac{I^n}{JI^{n-1}}\right) -
\lambda\left(\frac{I^{n-1} : J}{I^{n-2}}\right).
$$
\end{corollary}
\begin{proof}
Set $K = R$, $\F = \{I^n\}$ in Lemma \ref{funda}. Then $H_K(\F, n) =
H(I,n)$ for all $n \geq 0$ so that $P_K(\F,n) = P(I,n)$ and hence the
assertion follows.
\end{proof}

\noindent
As a consequence of the generalization of the Fundamental Lemma, we
obtain expressions for the Hilbert coefficients $g_1$ and $g_2$.
\begin{corollary}\label{formulaforgi}
Set 
$$
v_n =  \left\{
  \begin{array}{lll}
  e_0(\F) \mbox{ if } n = 0 \\
  e_0(\F) - \lambda(R/KI_1) + \lambda(R/K) \mbox{ if } n = 1 \\
  \lambda(KI_n/KJI_{n-1}) - \lambda(KI_{n-1} : J/KI_{n-2}) \mbox{ if }
  n \geq 2. \end{array}
 \right.
$$
Then $g_1 = \sum_{n\geq1}v_n$ and $g_2 = \sum_{n\geq1}(n-1)v_n + 
\lambda(R/K).$
\end{corollary}
\begin{proof}
From Lemma \ref{funda} we have, 
$$
\sum_{n\geq 0}\Delta^2[P_K(\F,n) - H_K(\F,n)]t^n
= \sum_{n\geq0}v_n t^n.
$$ Write $P_K(\F,n) = e_0(\F){n+2 \choose 2} -
g_1' (n+1) + g_2'$. Then comparing with the earlier notation, we get
$g_1 = g_1' - e_0(\F)$ and $g_2 = e_0(\F) - g_1' + g_2'$. Since
$P_K(\F,n)$ is a polynomial of degree 2, $\Delta^2P_K(\F,n) = e_0(\F)$
for all $n \geq 0$ so that $\sum_{n\geq0}\Delta^2P_K(\F,n)t^n =
e_0(\F)/(1-t).$ Let $\sum_{n \geq 0}H_K(\F,n) t^n = f(t)/(1-t)^3$.
Then by Proposition 4.1.9 of \cite{bh}, $e_0(\F) = f(1), ~~ g_1' = f'(1)$
and $g_2' = f''(1)/2!$. Also we have,
\begin{eqnarray*}
\sum_{n\geq0} \Delta^2 H_K(\F,n)t^n & = & \sum_{n\geq0} H_K(\F,n)t^n -
2 \sum_{n\geq0}H_K(\F,n-1) t^n + \sum_{n\geq0} H_K(\F,n-2)t^n \\
& = & \frac{f(t)}{(1-t)^3} - 2 H_K(\F, -1) - 2t\frac{f(t)}{(1-t)^3} +
H_K(\F,-2) + t H_K(\F,-1) + t^2 \frac{f(t)}{(1-t)^3} \\
& = & \frac{f(t)}{(1-t)} - 2 \lambda(R/K) + \lambda(R/K) +
t\lambda(R/K) \\
& = & \frac{f(t) - (1-t)^2\lambda(R/K)}{(1-t)}.
\end{eqnarray*}
Therefore 
$$
\sum_{n\geq0}\Delta^2 [P_K(\F,n) - H_K(\F,n)]t^n = \frac{e_0(\F) -
f(t) + (1-t)^2\lambda(R/K)}{(1-t)} = \sum_{n\geq0}v_nt^n.
$$
Therefore 
\begin{eqnarray*}
e_0(\F) - f(t) + (1-t)^2\lambda(R/K) = (1-t) \sum_{n\geq0} v_nt^n. 
\nonumber
\end{eqnarray*}
Hence
\begin{eqnarray}
f(t)  =   e_0(\F) + (1-t)^2\lambda(R/K) - (1-t) 
\sum_{n\geq0} v_nt^n.
\end{eqnarray}
Thus, $f(1) = e_0(\F)$. Differentiating (1) with respect to $t$, we 
get
$$
f'(t) = -2(1-t) \lambda(R/K) - (1-t)\sum_{n\geq0}nv_nt^{n-1} +
\sum_{n\geq0} v_nt^n.
$$
Therefore $g_1' = f'(1) = \sum_{n\geq0}v_n.$ Differentiating (1) twice
with respect to $t$, we get
$$
f''(t) = 2\lambda(R/K) - (1-t)\sum_{n\geq0}n(n-1)v_nt^{n-2} +
2\sum_{n\geq0}nv_nt^{n-1}
$$
so that 
$$
g_2' = f''(1)/2 = \sum_{n\geq0}nv_n + \lambda(R/K).
$$
Therefore 
$$
g_1 = g_1' - e_0(\F) = \sum_{n\geq0}v_n - e_0(\F) =
\sum_{n\geq1}v_n
$$
and 
$$
g_2 = g_2' - g_1 = \sum_{n\geq0}nv_n + 
\lambda(R/K) - \sum_{n\geq1}v_n = \sum_{n\geq1}(n-1)v_n + 
\lambda(R/K).
$$
\end{proof}

\noindent
{\bf Remark:} The above formulas for $g_1$ and $g_2$ generalize the 
formulas for $e_1$ and $e_2$ obtained by Huneke as consequences of 
his fundamental lemma. To obtain Huneke's formulas for $e_1$ and 
$e_2$, one simply puts $K=R$ in the above formulas for $g_1$ and 
$g_2$.

\begin{example}
{\em Let $k$ be any field and let $R = k[\![x,y]\!]$. Let $I = (x^3, x^2y,
y^3).$ Then $J = (x^3,  y^3)$ is a minimal reduction of $I$. Then
$r_J(I) = 3$, $\m I^n = \m JI^{n-1}$ for all $n \ge 2$ and $\m I^n :
J = \m I^{n-1}$ for all $n \geq 1$.
Therefore $v_n = 0$ for all $n \geq 2.$ One can also see that $e_0 =
9$ and $\lambda(R/\m I) = 10$. Hence we have
$v_0 = 9, v_1 = e_0 - \lambda(R/mI) + \lambda(R/m) = 9 - 10 + 1 = 0$
Thus $g_1' = v_0 + v_1 = 9$ and $g_2' = v_1 + \lambda(R/m) = 1$ which
gives $g_1 = g_1' - e_0 = 0$ and $g_2 = g_2' - g_1' + e_0 = 1$.}
\end{example}

\noindent
The following lemma shows that the behaviour of the superficial
elements in $F_K(\F)$  is quite similar to that of superficial elements 
in $G(\F)$. 

\begin{lemma}\label{modx}
Let $x$ be a regular element in $I_1$ such that $x^o$ is superficial
in $F_K(\F)$ and $x^*$ is superficial in $G(\F)$. Let $\bar{g_i}$
denote the coefficients of the polynomial corresponding to the
function $\lambda(\bar{R}/\bar{K}\bar{I_n})$, where $``-"$ denotes
$\modulo (x)$. Then $\bar{g_i} = g_i$ for all $i = 0, \ldots, d-1$.
\end{lemma}
\begin{proof}
Consider the exact sequence
$$
0 \longrightarrow \frac{KI_n : x}{KI_n} \longrightarrow R/KI_n
\stackrel{x}{\longrightarrow} R/KI_n \longrightarrow R/(KI_n + xR)
\longrightarrow 0.
$$
Then $\lambda (\bar{R}/\bar{K}\bar{I}_n) = \lambda(R/(KI_n+xR)) =
\lambda(KI_n : x/KI_n)$. Since $x^o$ is superficial in $F_K(\F)$ and
$x^*$ is superficial in $G(\F)$, $KI_n : x = KI_{n-1}$ for $n \gg 0$,
by Lemma \ref{supchar}. Hence,  $\lambda(R/(KI_n + xR)) =
\lambda(R/KI_n) - \lambda(R/KI_{n-1})$. for $n \gg 0$. Therefore 
\begin{eqnarray*}
P_{\bar{K}}(\bar{\F}, n) & = & P_K(\F, n) - P_K(\F, n-1) \\
& = & \sum_{i = 0}^d (-1)^i g_i {n+d-i-1 \choose d-i} - 
\sum_{i = 0}^d (-1)^i g_i {n+d-i-2 \choose d-i} \\
& = & \sum_{i = 0}^{d-1} (-1)^i g_i {n+d-i-2 \choose d-1-i}. 
\end{eqnarray*}
\end{proof}

\section{Cohen-Macaulay fiber cones}

In this section we obtain a lower bound for the Hilbert coefficient
$g_1$. We will characterize the Cohen-Macaulayness of $F_K(\F)$ 
in terms of $g_1.$ In this characterization, We need to assume
that $G(\F)$  has almost maximal depth. We will show by an example
that we need this assumption for such a characterization. Let 
$\gamma(\F)$ denote $\grade_{G(\F)_+}(G(\F))$. We begin by giving a
lower bound for the Hilbert coefficient $g_1$.

\begin{proposition}\label{bound}
Let $(R, \m)$ be a Cohen-Macaulay local ring of dimension d with
infinite residue field and let $J \subseteq I_1$ be a minimal reduction of
$I_1$ and let $K$ be an ideal such that $I_{n+1} \subseteq KI_n$. 
Then $g_1 \geq \sum_{n \geq 1}\lambda(KI_n + J/J) - \lambda(R/K)$. 
\end{proposition}

\begin{proof}
We apply induction on $d$. Let $d = 1$ and let $(x) = J$. 
For $i \geq 0$, from the exact sequence
$$
0 \longrightarrow \frac{(KI_{i+1} : x) \cap I_i}{KI_i} 
\longrightarrow I_i/KI_i \stackrel{x}{\longrightarrow} I_{i+1}/KI_{i+1}
\longrightarrow I_{i+1}/(KI_{i+1} + xI_i) \longrightarrow 0.
$$
it follows that 
$$
\lambda\left(\frac{I_{i+1}}{KI_{i+1}}\right) - \lambda\left(\frac{I_i}
{KI_i}\right) = \lambda\left(\frac{I_{i+1}}{(KI_{i+1} + xI_i)}\right) - 
\lambda\left(\frac{(KI_{i+1} : 
x) \cap I_i}{KI_i}\right) \hspace*{1cm}\cdots\cdots ~~~{\mbox (E_i) }.
$$
Summing up $E_0, E_1, \cdots, E_{n-1}$, we get
$$
\lambda(I_n/KI_n) - \lambda(R/K)  =  \sum_{i=1}^n
\lambda(I_i/(KI_{i} + xI_{i-1})) - \sum_{i=1}^{n}\lambda((KI_i : x)
\cap I_{i-1}/KI_{i-1}).
$$
Thus for all $n \gg 0$,
$$
f_0 = \lambda(R/K)  +  \sum_{i \geq 1}
\lambda(I_i/(KI_{i} + xI_{i-1})) - \sum_{i \geq 1}\lambda((KI_i : x)\cap 
I_{i-1}/KI_{i-1}).
$$
By Theorem 4.7 of \cite{hm}, 
$e_1(\F) = \sum_{i \geq 1} \lambda(I_i/xI_{i-1})$. Thus
\begin{eqnarray*}
g_1 & = & e_1(\F) - f_0 \\
& = & \sum_{i \geq 1} \lambda(I_i/xI_{i-1}) - \lambda(R/K)  -  
\sum_{i \geq 1} \lambda(I_i/(KI_{i} + xI_{i-1})) + 
\sum_{i \geq 1}\lambda((KI_i : x)\cap I_{i-1}/KI_{i-1})\\
& = & \sum_{i \geq 1} \lambda((KI_{i} + xI_{i-1})/xI_{i-1}) + \sum_{i
\geq 1} \lambda((KI_i : x)\cap I_{i-1}/KI_{i-1}) - \lambda(R/K).
\end{eqnarray*}
This implies that $g_1 \geq \sum_{i \geq 1} \lambda((KI_{i} +
xI_{i-1})/xI_{i-1}) - \lambda(R/K) \geq \sum_{i \geq 1}\lambda(KI_i
+ J/J) - \lambda(R/K).$
Hence the result is true for $d = 1$. Let us assume that $d > 1$ and
the assertion is true for $d-1$. Choose the
generators $x_1, \ldots, x_d$ of $J$ such that $x_1^o$ (resp. $x_1^* 
\in G(\F)$) is superficial in $F_K(\F)$ (resp. $G(\F)$). Let $`` - "$
denote images $\modulo (x_1)$. Then by Lemma \ref{modx}, $\bar{g_1} =
g_1$. By induction 
\begin{eqnarray*}
\bar{g_1} & \geq & \sum_{n \geq 1} \lambda(\bar{K}\bar{I_n} +
\bar{J}/\bar{J}) - \lambda(\bar{R}/\bar{K}) \\
& = & \lambda((KI_n + xR) + (J + xR)/(J + xR)) - \lambda(R/K) \\
& = & \sum_{n\geq 1}\lambda(KI_n + J/J) - \lambda(R/K).
\end{eqnarray*}
\end{proof}

As a consequence of the above bound, we derive a a bound for the
multiplicity of the fiber cone given in \cite[Theorem 2.1]{cpv}.

\begin{corollary}
Let $(R,\m)$ be a Cohen-Macaulay local ring of dimension $d > 0$ and infinite
residue field. Let $I$ be an $\m$-primary ideal. Then the multiplicity
of the fiber cone, $F(I)$ satisfies
$$
f_0 \leq e_1(I) - e_0(I) + \lambda(R/I) + \mu(I) - d + 1.
$$
\end{corollary}

\begin{proof}
Since $f_0 = e_1 - g_1$, by Proposition \ref{bound}, we get
\begin{eqnarray*}
f_0 & \leq & e_1 - \sum_{n\geq1}\lambda\left(\m I^n+J/J\right) + 1 \\
    & = & e_1 - \lambda(\m I/\m J) + 1 - 
          \sum_{n\geq2}\lambda\left(\m I^n+J/J\right) \\
    & = & e_1(I) - e_0(I) - d + \lambda(R/I) + \mu(I) + 1 - 
          \sum_{n\geq2}\lambda\left(\m I^n+J/J\right) \\
    & \leq & e_1(I) - e_0(I) + \lambda(R/I) + \mu(I) - d + 1.
\end{eqnarray*}
\end{proof}

Now we prove a characterization for Cohen-Macaulayness of the fiber
cone in terms of $g_1$.

\begin{theorem}\label{cm}
Let $(R,\m)$ be a Cohen-Macaulay local ring of dimension $d > 0$ with
$ R/\m $ infinite.  Let $\F = \{I_n\}$ be a Hilbert filtration
of $R$ and let $K$ be an ideal of $R$ containing $I_1$. Let $J$ be a
minimal reduction of $I_1$. Assume that $\gamma(\F) \geq d-1$. Then
$F_K(\F)$ is Cohen-Macaulay if and only if $g_1 = \sum_{n\geq
1}\lambda(KI_n + JI_{n-1}/JI_{n-1}) - \lambda(R/K)$.
\end{theorem}

\begin{proof}
Suppose that $F_K(\F)$ is Cohen-Macaulay. As $J+KI_1/KI_1$ is
generated by a homogeneous system of parameters of degree 1, 
$f_0 = \lambda(F_K(\F)/JF_K(\F)) = \sum_{i \geq 0}\lambda(I_i/(KI_i +
JI_{i-1}))$. Since $\gamma(\F) \geq d-1$, by Theorem 4.7 of \cite{hm}, 
$e_1 = \sum_{i \geq 1}\lambda(I_i/JI_{i-1})$. Therefore
\begin{eqnarray*}
g_1 & = & e_1 - f_0 = \sum_{n \geq 1}\lambda(I_n/JI_{n-1}) - 
\sum_{n \geq 0}\lambda(I_n/(KI_n + JI_{n-1})) \\
& = & \sum_{n \geq 1}\lambda(I_n/JI_{n-1}) - 
\sum_{n \geq 1}\lambda(I_n/(KI_n + JI_{n-1})) - \lambda(R/K) \\
& = & \sum_{n\geq 1}\lambda(KI_n + JI_{n-1}/JI_{n-1}) - 
\lambda(R/K).
\end{eqnarray*}
Conversely, suppose 
$g_1 = \sum_{n \geq 1}\lambda(KI_n + JI_{n-1}/JI_{n-1}) -
\lambda(R/K)$. 
Then by reversing the above steps, one can see that $f_0 = 
\sum_{n \geq 0}\lambda(I_n/(KI_n + JI_{n-1}))= \lambda(F_K(\F)/JF_K(\F)).$
 Therefore $F_K(\F)$ is Cohen-Macaulay.
\end{proof}
\vskip 2mm
\noindent
The following example shows that the assumption in Theorem
\ref{cm} that $\depth G(\F) \geq d-1$ cannot be dropped. 
\begin{example}{\em
Let $R = k[\![ x, y ]\!], \m = (x, y)$ and $I = (x^4, x^3y, xy^3,
y^4)$. Then $J = (x^4, y^4)$ is a minimal reduction of $I$ and $I^3 =
JI^2$. Note that $I^n = \m^{4n}$ for all $n \geq 2$. We compute the
Hilbert coefficients of $I$. Since $I^n = \m^{4n}$ for all $n \geq 2$,
$\lambda(R/I^n) = \lambda(R/\m^{4n}) = {4n + 1 \choose 2} = e_0(I){n+1
\choose 2} - e_1(I)~n + e_2(I)$. Solving the equation by putting
various values for $n$, we get $e_0(I) = 16,~ e_1(I) = 6,~ e_2(I) =
0$.  From direct computations one can see $\lambda(I/J) = 5$ and
$\lambda(I^2/JI) = 2$. Hence $e_1(I) <
\sum_{n\geq1}\lambda(I^n/JI^{n-1})$. Therefore $\depth G(I) = 0$.
Since $I^n = \m^{4n}$ for all $n \geq 2$, $\mu(I^n) = \lambda(\m^{4n}
/\m^{4n+1}) = {4n+1 \choose 1}$ for all $n \geq 1$. Therefore $f_0 =
4$ so that $g_1 = e_1 - f_0 = 2$. Also, one can see that $\lambda(\m
I+J/J) \neq 0$ and $\m I^n \subseteq JI^{n-1}$ for all $n \geq 2$.
Then $\lambda(\m I+J/J) = \lambda(\m I/\m I \cap J) = \lambda(\m I/\m
J) = 3.$ Therefore $\sum_{n\geq1} \lambda(\m I^n + JI^{n-1}/JI^{n-1})
- 1 = 2 = g_1.$ The Hilbert Series of the fiber cone is given by 
$$
H(F(I), t) = \frac{1+2t+2t^2-t^3}{(1-t)^2}.
$$ 
Since the numerator contains a negative coefficient, $F(I)$ is not
Cohen-Macaulay. 

}
\end{example}

\section{Fiber cones with almost maximal depth}
In this section we present a characterization for the fiber cone $F_K(\F)$
to have almost maximal depth in terms of $g_1.$ This is an analogue of the 
Huckaba-Marley characterization for the associated graded ring to have 
almost maximal depth referred above. We will need to invoke a result
due to Cortadellas and Zarzuela from \cite{cz} which gives a criterion 
for a sequence of degree one elements in  $F_K(\F)$ to be a regular sequence.
We take this opportunity to present a simple proof of their result since it
is a very basic result for detecting regular sequences in fiber cones.

\begin{theorem}\label{fiberdepth}
Let $(R,\m)$ be a Noetherian local ring, $\F$ a filtration of ideals,
$K$ an ideal containing $I_1$ and $x_1, \ldots, x_k \in I_1$. Assume 
that 
\begin{enumerate}
\item[(i)]  $x_1, \ldots, x_k$ is a regular sequence in $R$.
\item[(ii)] $x_1^* \ldots, x_{k}^* \in G(\F)$ is a regular sequence.
\item[(iii)] $x_1^o, \ldots, x_k^o \in F_K(\F)$ is a superficial 
sequence.
\end{enumerate}
Then $\depth_{(x_1^o, \ldots, x_k^o)}F_K(\F) = k$ if and only if
$(x_1, \ldots, x_k) \cap KI_n = (x_1, \ldots, x_k)KI_{n-1}$ for
all $n \geq 1$.
\end{theorem}

\begin{proof}
We induct on $k$. Let $k = 1$. Let $(x) \cap K I_{n} = xK I_{n-1}$ for
all $n \geq 1$. Then $K I_{n} : x = K I_{n-1}$ for all $n \geq 1$ and
hence $x^o$ is regular in $F_K(\F)$. Suppose $x^o$ is regular in $F_K(\F)$.
Let $n \geq 1$ and $y \in (x) \cap K I_n$. Let $y = rx$ for some $r
\in R$. Then $r \in K I_n : x \subseteq I_n : x$. Since $x^*$ is a
nonzerodivisor in $G(\F)$, $I_n : x = I_{n-1}$ for all $n \geq 1$.
Therefore $r \in I_{n-1}$. Hence $r \in (K I_n : x) \cap I_{n-1} = K
I_{n-1}$ so that $y \in xK I_{n-1}$. 
Now assume that $k > 1$ and the result is true for all $l \leq k-1$.
Put $J = (x_1, \ldots, x_k)$, $J^o = (x_1^o, \ldots, x_k^o) \subseteq
F_K(\F)$ and $J^* = (x_1^*, \ldots, x_k^*) \subseteq G(\F)$.  Let
$``-"$ denote images modulo ($x_1$). Then $F_K(\F)/x_1^oF_K(\F)$ $\cong
F_{\bar{K}}(\bar{\F})$.  Assume $J \cap K I_n = JKI_{n-1}$ for all $n \geq 1$.
Then 
$$
\bar{J} \cap \bar{K}\bar{I}_n = (J + x_1R) \cap (K I_n + x_1R) =
J \cap (K I_n + x_1R) = (J \cap K I_n) + x_1R = JK I_{n-1} + x_1R =
\bar{J}\bar{K}\bar{I}_{n-1}.
$$
By induction
$\depth_{\bar{J}^o}F_{\bar{K}}(\bar{\F}) = k-1 $. Thus $x_1^o$ is regular in
$F_K(\F)$ and hence, 
$\depth_{J^o}(F_K(\F)) = k$, by Lemma \ref{sally}. Conversely assume
that $\depth_{J^o}(F_K(\F)) = k$. Since $x_1^o$ is superficial in
$F_K(\F)$ and $\depth F_K(\F) > 0$, $x_1^o$ is regular in $F_K(\F)$.
Then $\depth_{\bar{J}^o}(F(\bar{I})) = k-1$. Hence, by induction,
$\bar{J} \cap \bar{K}\bar{I}_n = \bar{J}\bar{K}\bar{I}_{n-1}$.
Therefore, $J \cap K I_n + x_1R = JK I_{n-1} + x_1R$. Hence
\begin{eqnarray*}
  J \cap K I_n & = & J K I_{n-1} + (x_1R \cap (J \cap K I_n)) \\
    & = & J K I_{n-1} + (x_1R \cap K I_n)  \\
   & = & J K I_{n-1} + (x_1K I_{n-1}) 
    =  J K I_{n-1}.
\end{eqnarray*}
Therefore $J \cap K I_n = J K I_{n-1}$.
\end{proof}
\vskip 2mm
\noindent
We need the following lemma in the proof of the characterization for
the fiber cone to have depth at least $d-1$.

\begin{lemma}\label{regseq}
Let $(R,\m)$ be a $d$-dimensional Cohen-Macaulay local ring and let
$\F = \{I_n\}$ be a Hilbert filtration of $R$ such that $\gamma(\F)
\geq d-1$. Let $K$ be an ideal of $R$ containing $I_1$. Let $J = (x_1,
\ldots, x_d)$ be a minimal reduction of $I_1$ such that $x_1^*, \ldots,
x_{d-1}^*$ is a regular sequence in $G(\F)$.  If $KI_n \cap (x_1,
\ldots, x_{d-1}) \subseteq JKI_{n-1}$ for all $n \geq 1$, then $x_1^o,
\ldots, x_{d-1}^o$ is a regular sequence in $F_K(\F)$.
\end{lemma}
\begin{proof}
Since $x_1^*, \ldots, x_{d-1}^*$ is a regular sequence in $G(\F)$, by
Theorem \ref{fiberdepth}, it is enough to show that $KI_n \cap (x_1,
\ldots, x_{d-1}) = (x_1, \ldots, x_{d-1})KI_{n-1}$ for all $n \geq 1$.
Induct on $n$.  Let $z \in KI_1 \cap (x_1, \ldots, x_{d-1})$.  Write
$z = \sum_{i=1}^{d-1}r_ix_i$ for $r_i \in R$.  Since $KI_1 \cap (x_1,
\ldots, x_{d-1}) \subseteq (x_1, \ldots,
x_{d})K$, $z = \sum_{i=1}^{d-1}s_ix_i + px_d$ for some $p, s_i \in K$.
Then $px_d \in (x_1, \ldots, x_{d-1})$ and hence $p \in (x_1, \ldots,
x_{d-1})$. Since $x_d \in K$, $z \in
(x_1, \ldots, x_{d-1})K$. Now assume that $n \geq 2$ and for all $l <
n$, 
$$
KI_l \cap (x_1, \ldots, x_{d-1}) = (x_1, \ldots, x_{d-1})KI_{l-1}.
$$ 
Let $z \in KI_n \cap (x_1, \ldots, x_{d-1}) \subseteq JKI_{n-1}$.
Write $z = \sum_{i=1}^{d-1} r_ix_i = \sum_{i=1}^{d-1}s_ix_i + px_d$,
where $r_i \in R, ~~p, s_i \in KI_{n-1}$. Then $px_d \in (x_1, \ldots,
x_{d-1})$ and hence $p \in (x_1, \ldots, x_{d-1})$. Therefore 
$$
KI_{n} \cap (x_1, \ldots, x_{d-1}) = (x_1, \ldots, x_{d-1})KI_{n-1} +
x_d(KI_{n-1} \cap (x_1, \ldots, x_{d-1})).
$$
By induction $KI_{n-1}
\cap (x_1, \ldots, x_{d-1}) \subseteq (x_1, \ldots, x_{d-1})KI_{n-2}$.
Hence $p \in (x_1, \ldots, x_{d-1})KI_{n-2}$ so that 
$px_d \in (x_1,\ldots, x_{d-1})KI_{n-1}$.  Therefore 
$z \in (x_1, \ldots, x_{d-1})KI_{n-1}$.
\end{proof}

\noindent
We prove a necessary and sufficient condition for the fiber cone to
have depth at least $d-1$ in terms of $g_1$ in the following theorem.
\begin{theorem}\label{fiberd-1}
Let $(R,\m)$ be a Cohen-Macaulay local ring of dimension $d > 0$. Let
$\F$ be a Hilbert filtration, $K$ an ideal such that $I_{n+1}
\subseteq KI_n$ for all $n \geq 0$ and $J$ a minimal reduction of
$I_1$. Assume that $\gamma(\F) \geq d-1$. Then $g_1 = \sum_{n\geq1}
\lambda(KI_n/JKI_{n-1}) - \lambda(R/K)$ if and only if $\depth F_K(\F)
\geq d-1$. 
\end{theorem}
\begin{proof}
We induct on $d$. Let $d = 1$. In this case, we need to show that
$g_1$ is equal to the summation given in the statement of the theorem.
Let $(x)$ be a reduction of $I_1$. From
the proof of Proposition \ref{bound} we get, 
$$
g_1 = \sum_{n \geq 1} \lambda((KI_{n} + xI_{n-1})/xI_{n-1}) + 
\sum_{n \geq 1} \lambda((KI_n : x) \cap I_{n-1}/KI_{n-1}) - 
\lambda(R/K).
$$
{\it Claim :} For all $n \geq 1$, $(KI_n : x) \cap I_{n-1}/KI_{n-1}
\cong xI_{n-1} \cap KI_n /xKI_{n-1}$.\\
Consider the multiplication map $\mu_x : (KI_n : x) \cap 
I_{n-1}/KI_{n-1} \longrightarrow xI_{n-1} \cap KI_n /xKI_{n-1}$.
Let $y = xs \in xI_{n-1} \cap KI_n$ for some $s \in I_{n-1}$. Then $s
\in (KI_n : x) \cap I_{n-1}$. Therefore, $\mu_x$ is surjective.
Let  $y \in (KI_n : x) \cap I_{n-1}$ and  $xy \in xKI_{n-1}$.
Since $x$ is regular in $R$, $y \in KI_{n-1}$. Hence $\mu_x$ is
injective so that $\mu_x$ is an isomorphism. Therefore we have,
\begin{eqnarray*}
g_1 & = & \sum_{n\geq1} [\lambda(KI_n/(KI_n \cap 
xI_{n-1})) + \lambda((xI_{n-1} \cap KI_n)/xKI_{n-1})] - \lambda(R/K) \\
& = & \sum_{n\geq1} \lambda(KI_n/xKI_{n-1}) - \lambda(R/K).
\end{eqnarray*}
Assume now that $d > 1$. Choose $x_1, \ldots, x_d$ such that $J =
(x_1, \ldots, x_d),\; x_1^*, \ldots, x_d^*$ is a superficial sequence in
$G(\F)$ and $(x_1^o, \ldots, x_d^o)$ is a superficial sequence in
$F_K(\F)$. Since $\gamma(\F) \geq d-1$, $x_1^*, \ldots, x_{d-1}^*$ is
a regular sequence in $G(\F)$ (existence of such a generating set can
be derived from Proposition A.2.4 of \cite{ma}).  Suppose $g_1 =
\sum_{n\geq1} \lambda(KI_n/JKI_{n-1}) - \lambda(R/K)$.  Let $``-"$
denote images modulo $(x_1, \ldots, x_{d-1})$. Then by Lemma
\ref{modx}, $g_1 = \bar{g_1}$ and  
\begin{eqnarray*}
\bar{g_1} & = & \sum_{n\geq1} 
\lambda(\bar{K}\bar{I}_n/\bar{J}\bar{K}\bar{I}_{n-1}) -
\lambda(\bar{R}/\bar{K}) \hspace*{1cm} \mbox{(Since $\dim \bar{R}$ =
1)}\\
& = & \sum_{n\geq1} \lambda((KI_n + (x_1, \ldots,
x_{d-1}))/(JKI_{n-1} + (x_1, \ldots,x_{d-1}))) - \lambda(R/K) \\
& = & \sum_{n\geq1} \lambda(KI_n/(JKI_{n-1} + KI_n\cap(x_1, \ldots, 
x_{d-1}))) - \lambda(R/K). 
\end{eqnarray*}
By assumption $g_1 = \sum_{n\geq1}\lambda(KI_n/JKI_{n-1}) -
\lambda(R/K)$. Therefore, $(x_1, \ldots, x_{d-1}) \cap KI_n \subseteq
JKI_{n-1}$. Hence by Lemma \ref{regseq}, $x_1^o, \ldots, x_{d-1}^o$ is
a regular sequence in $F_K(\F)$.
\vskip 2mm
\noindent
Conversely, let $\depth F_K(\F) \geq d-1$. Choose $x_1 \in I_1$ such
that $x_1^*$ is regular in $G(\F)$ and $x_1^o$ is regular in
$F_K(\F)$. Hence $F_K(\F)/x_1^oF_K(\F) \cong F_{\bar K}(\bar{\F})$
and $g_1 = \bar{g_1}$, where $`` - "$ denote images modulo $(x_1)$.
Then $\depth F_K(\F)/x_1^oF_K(\F) \geq d-2$. By induction 
\begin{eqnarray*}
\bar{g_1} & = & \sum_{n\geq1}
\lambda(KI_n + x_1R/JKI_{n-1}+x_1R) - \lambda(R/K) \\
& = & \sum_{n\geq1}
    \lambda(KI_n/(JKI_{n-1} + (x_1R \cap KI_{n}))) - \lambda(R/K). 
\end{eqnarray*}

Since $x_1^o$ is regular in $F_K(\F)$ and $x_1^*$ is regular in $G(\F)$, 
$x_1R \cap KI_n = x_1KI_{n-1}$. Therefore
$$
\bar{g_1} = \sum_{n\geq1} \lambda(KI_n/(JKI_{n-1} + x_1KI_{n-1}))
- \lambda(R/K) = 
\sum_{n\geq1} \lambda(KI_n/JKI_{n-1}) - \lambda(R/K).
$$
Since $g_1 = \bar{g_1}$, the assertion follows.
\end{proof}

\begin{remark}\label{upperbound}{\em
From the above proof we obtain an upper bound for the coefficient
$g_1$:
$$
g_1 \leq \lambda\left(KI_n/KJI_{n-1}\right) - \lambda(R/K).
$$
}
\end{remark}

\section{Cohen-Macaulay fiber cones of ideals with minimal 
multiplicity}

Let $(R,\m)$ be a Cohen-Macaulay local ring of dimension $d > 0$. Let
$I$ be an $\m$-primary ideal of $R$ and $J$ be a minimal reduction of
$I$. Let $K$ be an ideal containing $I$. Let $F_{K}(I)$ be the fiber
cone of $I$ with respect to $K$ and let $G(I)$ be the associated
graded ring of $I$. Let $\mu(I)$ denote the minimum number of
generators of $I$. It is known that $e_0(I) + d - \lambda(R/I) \geq
\mu(I)$, \cite{g}. Goto defined an ideal $I$ to have minimal
multiplicity if $e_0(I) + d - \lambda(R/I) = \mu(I)$. It can be seen
that $I$ has minimal multiplicity if and only if for any minimal
reduction $J$ of $I$, $I\m = J\m$.   We generalize this notion. An ideal
$I$ is said to have minimal multiplicity with respect to $K$ if $KI =
KJ$ for any minimal reduction $J$ of $I$. Let $H_{K}(I, n) = 
\lambda(R/KI^n)$ and 
$P_{K}(I,n) = \sum_{i=0}^d (-1)^i g_i {n+d-i-1 \choose d-i}$ 
be the corresponding polynomial.

\begin{proposition}\label{mmg1}
Let $(R, \m)$ be a Cohen-Macaulay local ring of dimension $d$.\\
(i) If $I$ has minimal multiplicity with respect to $K$, 
then $g_1 = -\lambda(R/K)$. \\
(ii) If $KI\cap J =KJ$ and $g_1=-\lambda(R/K)$ then 
     $I$ has minimal multiplicity. \\
(iii) The ideal $I$ has minimal multiplicity if and only if $g_1=-1.$
  
\end{proposition}
\begin{proof}
(i) Suppose $I$ has minimal multiplicity with respect to $K$. Let $J$ be a
minimal reduction of $I$. Then $KI^n = KJ^n$ for all $n \geq 1$.
Therefore for all $n \geq 1$,
\begin{eqnarray*}
\lambda(R/K I^n) & = & \lambda(R/K J^n) = \lambda(R/J^n) +
\lambda(J^n/K J^n) \\
& = & e_0(I){n+d-1 \choose d} + \lambda(R/K) {n+d-1 \choose d-1}.
\end{eqnarray*}
Hence $g_1 = -\lambda(R/K)$. \\
(ii) Suppose $g_1=-\lambda(R/K)$ and $KI\cap J=KJ.$ 
By Proposition \ref{bound},
$g_1 \geq \sum_{n \geq 1}\lambda(KI^n+J/J)-\lambda(R/K)$. Hence
$KI\subseteq J.$ Thus $KI=KI\cap J=KJ.$ Hence $I$ has minimal multiplicity 
with respect to $K.$ \\
(iii) Follows from (i) and (ii) since $\m I \cap J=\m J.$
\end{proof}

In the following Proposition we generalize, the equivalence of 
Cohen-Macaulayness of the fiber cone and the almost maximal depth for
$G(I)$, given in Proposition 2.5 of \cite{g}.

\begin{proposition}\label{mmcmfiber}
Let $(R,\m)$ be a Cohen-Macaulay local ring, $I$ an $\m$-primary ideal, $J$ a
minimal reduction of $I$ and $K$ an ideal containing $I$. Suppose that
there exists an integer $t$ such that
$$
KI^n \cap J = KJI^{n-1} \mbox{ for } n = 1, \ldots, t \;\;\; 
\mbox{ and }\;\;\; KI^{t+1} = KJI^t. \hspace*{1in} (*)
$$
Then $F_K(I)$ is Cohen-Macaulay if and only if $\gamma(I) \geq d-1$.
\end{proposition}
\begin{proof}
From $(*)$, it follows that $KI^n \cap J = KI^n \cap JI^{n-1} =
KJI^{n-1}$ for all $n = 1, \ldots, t$. Therefore, from Remark
\ref{upperbound} and Proposition \ref{bound}, it follows that 
$$
g_1 = \sum_{n\geq1}\lambda\left(KI^n/KJI^{n-1}\right) - \lambda(R/K)
= \sum_{n\geq1}\lambda\left(KI^n + JI^{n-1}/JI^{n-1}\right) - \lambda(R/K).
$$
Therefore by Theorem \ref{cm}, if $\gamma(I) \geq d-1$, then $F_K(I)$
is Cohen-Macaulay.
\vskip 2mm
\noindent
Assume that $F_K(I)$ is Cohen-Macaulay. Then $J+KI/KI$ is system of
parameters for $F_K(I)$ and 
$$
f_0 = \lambda(F_K(I)/JF_K(I)) = \sum_{n\geq0}\lambda(I^n/KI^n+JI^{n-1})
    = \lambda(R/K) + \sum_{n\geq1}\lambda(I^n/KI^n+JI^{n-1}).
$$
Hence
\begin{eqnarray*}
e_1(I) & = & f_0 + g_1 \\
       & = & \lambda(R/K) +
             \sum_{n\geq1}\lambda(I^n/KI^n+JI^{n-1}) + \sum_{n\geq1}
             \lambda\left(KI^n + JI^{n-1}/JI^{n-1}\right) - 
	     \lambda(R/K) \\
       & = & \sum_{n\geq1}\lambda(I^n/JI^{n-1}).
\end{eqnarray*}
Therefore, by Theorem 4.7 of \cite{hm}, $\gamma(I) \geq d-1$.
\end{proof}

\begin{corollary}\label{cmd-1}
Let $(R,\m)$ be a $d$-dimensional Cohen-Macaulay local ring and let 
$I$ be an $\m$-primary ideal with minimal multiplicity with respect to
$K \supseteq I$. Suppose $KI\cap J = KJ$ Then $F_K(I)$ is
Cohen-Macaulay if and only if $\gamma(I) \geq d-1$.
\end{corollary}

\begin{theorem}\label{gotomm}
Let $(R,\m)$ be a Cohen-Macaulay local ring of dimension $d$ and let
$I$ be an $\m$-primary ideal of $R$ with minimal multiplicity. Then
the following statements are equivalent:
\begin{enumerate}
\item[$1.$] $G(I)$ is Cohen-Macaulay.
\item[$2.$] $F(I)$ is Cohen-Macaulay and $r(I) \leq 1.$
\item[$3.$] $r(I) \leq 1$.
\end{enumerate}
\end{theorem}
\begin{proof}
(1) $\Rightarrow$ (2) : Since $I$ has minimal multiplicity, $I\m = 
J\m$ for any minimal reduction $J$ of $I$. Assume that $G(I)$ is
Cohen-Macaulay. Then by Proposition \ref{cmd-1}, $F(I)$ is
Cohen-Macaulay. Therefore 
\begin{eqnarray*}
f_0 & = & \sum_{n\geq 0}\lambda(I^n /JI^{n-1} + \m I^n) = 1 + 
\sum_{n \geq 1} \lambda(I^n/JI^{n-1}) \mbox{ (since } \m I^n \subseteq
JI^{n-1} \mbox{ for all } n \geq 1)\\
& = & 1 + \lambda(I/J) + \sum_{n\geq2} \lambda(I^n/JI^{n-1}) \\
& = & 1 + e_0(I) - \lambda(R/I) + \sum_{n\geq2} \lambda(I^n/JI^{n-1})
\\
& = & 1 + \mu(I)  - d + \sum_{n\geq2} \lambda(I^n/JI^{n-1})
\end{eqnarray*}
Since $I^2 \subseteq I\m = J\m \subseteq J$, $I^2 = I^2 \cap
J = IJ$, as $G(I)$ is Cohen-Macaulay. 

\noindent
(2) $\Rightarrow$ (3) : Clear.

\noindent
(3) $\Rightarrow$ (1) : This is known but we recall the argument.
Assume that $r(I) \leq 1$. Then $I^{n+1} = JI^n$ for all $n \geq 1$ so
that $JI^n = I^{n+1} \cap J$ for all $n \geq 1$. Therefore by
\cite{vv}, $x_1^*, \ldots, x_d^*$ is a regular sequence in $G(I)$. 
Hence $G(I)$ is Cohen-Macaulay.
\end{proof}

\begin{example}{\em 
Let $R = k[\![t^4, t^5, t^6, t^7]\!]$ and let $I = (t^4, t^5, t^6).$
Then $J = (t^4)$ is a minimal reduction of $I$. It can easily be
checked that $I\m = J\m.$ Hence $F(I)$ is Cohen-Macaulay. 
Since $I$ has
minimal multiplicity, $g_1 = -1$. Since $\m I^n \subseteq J$ for all
$n \geq 1$, $g_1 = -1 = \sum_{n\geq1}\lambda(\m I^n +
JI^{n-1}/JI^{n-1}) - \lambda(R/\m)$.
}
\end{example}


\begin{thebibliography}{AAAA}

\bibitem [BH]{bh} W. Bruns and J. Herzog, {\em Cohen-Macaulay
rings}, Revised Edition, Cambridge Studies in Advanced Mathematics, 
39. Cambridge University Press, Cambridge, 1998
 
\bibitem[C]{c} A. Corso, {\em Sally modules of m-primary ideals in
  local rings}, preprint, arXiv:math.AC/0309027. 

\bibitem[CGPU]{cgpu} A. Corso, L. Ghezzi, C. Polini and B. Ulrich,
        {\em Cohen-Macaulayness of special fiber rings }Comm. Algebra 31 
        (special issue in honor of S. Kleiman) (2003), 3713-3734. 

\bibitem[CPV]{cpv} A. Corso, C. Polini and W. V. Vasconcelos, 
{\em Multiplicity of the special fiber ring of blowups}, preprint,
arXiv:math.AC/0307037

\bibitem [CZ]{cz} T. Cortadellas and S. Zarzuela, {\em On the
depth of the fiber cone of filtrations},  J. Algebra  {\bf 198} 
(1997),  no. 2, 428--445.

\bibitem [DRV]{drv} C. D'Cruz, K. N. Raghavan and J. K. Verma, 
{\em Cohen-Macaulay fiber cones}, Commutative Algebra, Algebraic
Geometry and Computational Methods (Hanoi, 1996), 233--246, Springer, 
Singapore, 1999

\bibitem [DV]{dv} C. D'Cruz and J. K. Verma, {\em Hilbert series of
  fiber cones of almost minimal mixed multiplicity}, J. Algebra {\bf
  251} (2002), no. 1, 98--109.

\bibitem [G]{g} S. Goto, {\em Cohen-Macaulayness and negativity of
  $A$-invariants in Rees algebras associated to $\m$-primary ideals
  of minimal multiplicity}, Commutative algebra, homological algebra
  and representation theory (Catania/Genoa/Rome, 1998).  J. Pure Appl.
  Algebra  {\bf 152}  (2000),  no. 1-3, 93--107.

\bibitem[HM]{hm} S. Huckaba and T. Marley, {\em Hilbert coefficients and
depth of associated graded rings}, J. London Math. Soc. (2) {\bf 56}
(1997), 64-76.

\bibitem[H]{h2} C. Huneke, {\em Hilbert functions and symbolic powers}, 
Michigan Math. J. {\bf 34} (1987), no. 2, 293--318. 

\bibitem[HH]{hh} R. H\"ubl and C. Huneke,
{ \em Fiber cones and the integral closure of ideals}, 
Collect. Math. 52 (2001), no. 1, 85--100.


\bibitem[Ma]{ma} T. Marley, {\em Hilbert functions of ideals in 
Cohen-Macaulay rings}, Thesis, Purdue University (1989).
 
\bibitem [S1]{s1} J. D. Sally, {\em On the associated graded ring of 
a local Cohen-Macaulay ring}, J. Math. Kyoto Univ.,{\bf 17}(1977), 19-21.

\bibitem [S2]{s5} J. D. Sally, {\em Cohen-Macaulay local rings of
embedding dimension $e+d-2$}, J. Algebra {\bf 83} (1983) 393-408.

\bibitem [VV]{vv} P. Valabrega and G. Valla, {\em Form rings 
and regular sequences},  Nagoya Math. J.  {\bf 72} (1978), 93--101.
\end{thebibliography}
\end{document}